\documentclass[amsa]{ipart}

\RequirePackage{hyperref,algorithmicx, algorithm}
\usepackage{hyperref,bm,upgreek,enumitem}
\usepackage{amsmath, amsthm, amssymb,  amsfonts, amssymb, graphicx }

\usepackage{multirow}%
\usepackage{booktabs}

\usepackage{easyReview}
\usepackage{enumitem}
\usepackage{caption}
\usepackage{subcaption}
\usepackage{amscd}
\usepackage{mathrsfs}%
\usepackage[title]{appendix}%
\usepackage{xcolor}%
\usepackage{textcomp}%
\usepackage{manyfoot}%
\usepackage{algcompatible}
\usepackage{algpseudocode}%
\usepackage{listings}%
\usepackage{array}
\usepackage{float} 
\usepackage{amscd}
\usepackage{todonotes}
\allowdisplaybreaks

\startlocaldefs
\theoremstyle{plain}

\endlocaldefs
\newtheorem{Theorem}{Theorem}[section]

\newtheorem{remark}{Remark}%

\newtheorem{definition}{Definition}%

\newtheorem{note}{Note}
\numberwithin{equation}{section}

\newcommand{\Idx}[1]{\mathscr{K}_{#1}}        
\newcommand{\Rplus}[1]{\mathbb{R}_{+}^{#1}}   
\newcommand{\Rpp}[1]{\mathbb{R}_{++}^{#1}}    
\newcommand{\xlm}{z^*}                    

\newcommand{\obj}[1]{\psi_{#1}}               
\newcommand{\fvec}{\boldsymbol{\psi}}         

\newcommand{\Feas}{\mathcal{Z}}               

\newcommand{\mult}{\lambda}              
\newcommand{\muc}{\mu}                       

\newcommand{\tun}{\boldsymbol{\mathbb{T}}}                 
\newcommand{\tunobj}[1]{\mathbb{T}_{#1}}               
\begin{document}

\begin{frontmatter}
\title[A Tunneling Method for Nonlinear MOP]{A Tunneling Method for Nonlinear Multi-objective Optimization Problems}
\begin{aug}
    \author{\fnms{Bikram} \snm{Adhikary}\ead[label=e1]{bikram.adhikary333@gmail.com}} \and
    \author{\fnms{Md Abu Talhamainuddin} \snm{Ansary} \thanksref{t2} \ead[label=e2]{md.abutalha2009@gmail.com}}
    \address{Department of Mathematics\\Indian Institute of Technology Jodhpur\\ Jodhpur, India-342030\\
            B. Adhikary \printead{e1}\\M. A. T. Ansary \printead{e2}}
    \thankstext{t2}{Corresponding author}
\end{aug}
%

\begin{abstract}
In this paper, a tunneling method is developed for nonlinear multi-objective optimization problems using some ideas of the single objective tunneling method. The proposed method does not require any a priori chosen parameters or ordering information of the objective functions. At any critical point, an auxiliary function is developed to find a different critical point that dominates the previous one. By repeatedly applying the tunneling procedure, it is possible to construct a broader approximation to the global Pareto front in nonconvex multi-objective optimization problems that may contain multiple local Pareto fronts. An algorithm is then designed based on this auxiliary function, and the convergence of this algorithm is justified under some mild assumptions. Finally, several numerical examples are presented to illustrate the effectiveness of the proposed method and to justify the theoretical results.
\end{abstract}



\end{frontmatter}
\section{Introduction}\label{sec1}
Multi-objective optimization seeks the simultaneous minimization of multiple conflicting objectives and has become an essential framework in diverse fields such as environmental modeling, space mission design, portfolio management, decision sciences, healthcare applications etc.. Classical solution strategies are typically divided into scalarization-based techniques~\cite{ghane2015new,kasimbeyli2013conic,kasimbeyli2019comparison,Miettinen2012-rv} and heuristic approaches~\cite{castro2018evaluating,konak2006multi,mostaghim2007multi,laumanns2002combining}. Despite their widespread use, scalarization methods remain sensitive to user-defined parameters, whereas heuristic approaches lack guarantees of monotonic descent and thus cannot ensure convergence to a Pareto-optimal solution.
\par Recent advancements in multi-objective optimization have introduced gradient descent methods that extend single objective descent techniques to multi-objective problems. These methods have been applied to smooth unconstrained problems \cite{ansary2015modified, fliege2009newton, fliege2000steepest}, constrained optimization problems \cite{ansary2020sequential, ansary2021sqcqp, fliege2016sqp}, non-smooth multi-objective optimization problems \cite{ansary2023proximal, ansary2022proximal, bento2014proximal, tanabe2019proximal, tung2025second, upadhyay2025inexact} and uncertain optimization problems \cite{bai2024convergence, ghosh2023newton, kumar2024steepest, kumar2023newton}. While these approaches provide global or local convergence justifications under specific theoretical assumptions, they often struggle when applied to non-convex problems, where local efficient solutions dominate.

A key challenge in multi-objective optimization is constructing a well-distributed Pareto front, particularly in non-convex scenarios. Existing multi-start and spreading techniques \cite{ansary2021sqcqp, fliege2016sqp} attempt to improve Pareto front diversity, yet they frequently converge to locally efficient solutions rather than global ones. In non-convex multi-objective problems, conventional line-search techniques often fail to generate a global Pareto front unless the initial approximation is very close to a global efficient solution.
\par Global optimization techniques for non-convex single objective problems \linebreak
~\cite{chen2024new, ketfi2014global, tunneling1985, duanli2010global, pandiya2021non, renpu1990filled, yang2006new, yilmaz2019new} have provided important conceptual foundations for overcoming such difficulties. Building on these principles, this study introduces a novel multi-objective tunneling algorithm specifically designed for complex, non-convex multi-objective optimization problems. The proposed method enables structured transitions from local to global Pareto-optimal solutions by employing a multi-objective variant of an auxiliary function, without requiring any priori chosen parameter or ordering information of the objective functions. The key contributions of this research include:
\begin{itemize}
    \item introduction of a {multi-objective auxiliary function} designed to enhance solution exploration in complex, non-convex optimization problems;
    \item development of a {multi-objective tunneling method} that enables efficient transitions from local to global Pareto-optimal solutions;
    \item establishment of a structured global search mechanism, where the algorithm consistently escapes local Pareto-optimal solutions while refining promising search regions, ensuring an effective balance between exploration and exploitation;
    \item theoretical validation together with extensive empirical evaluation, demonstrating the superior performance of the proposed method across diverse multi-objective problem settings.
\end{itemize}
\par The subsequent sections of this paper are organized as follows.  Section \ref{sec_pre} describes the essential preliminaries and theoretical foundations relevant to the proposed methodology.  Section \ref{sec_tun} presents the multi-objective tunneling approach, which integrates a novel multi-objective auxiliary function to enhance global optimization efficacy.  This section proposes an algorithm and justifies its convergence.  Section \ref{sec_exe} presents comprehensive numerical experiments undertaken to demonstrate the efficacy of the suggested method. Finally, Section~\ref{(sec_con)} provides concluding remarks and discusses potential directions for future research.
\section{Preliminaries}
\label{sec_pre}
Throughout the paper, we use the notations $\Idx{p} := \{1, 2, \dots, p\}$,\\
$\Rplus{p} := \{x \in \mathbb{R}^p \mid x_i \geq 0,~ \forall i \in \Idx{p}\}$, 
$\Rpp{p} := \operatorname{int}(\Rplus{p})$ where $p \in \mathbb{N}$. 
For $x^1,x^2 \in \mathbb{R}^p$, vector inequalities are interpreted component-wise.

Let us consider a box constrained multi-objective optimization problem:
\begin{flalign*}
(MOP_{BC}): \underset{z \in\mathbb{R}^n}{\min}&\quad \fvec(z) = \left( \obj{1}(z), \obj{2}(z), \dots, \obj{m}(z) \right)\\
\text{s.t.}& \quad lb \leq z \leq ub,
\end{flalign*}
where \( m \geq 2 \) and \( \obj{j}: \mathbb{R}^{n} \to \mathbb{R},~ j \in \Idx{m} \) are assumed to be continuously differentiable. The set of feasible solutions of $\mathrm(MOP_{BC})$ is \linebreak
$\Feas := \left\{z \in \mathbb{R}^n ~|~ lb \leq z \leq ub\right\}.$ If any \( z^* \in \Feas \) simultaneously minimizes all objective functions, then $z^*$ is referred as an ideal solution. Such a solution, however, is rarely attainable in practice, since objectives typically conflict with one another: improving the performance of one objective often comes at the expense of another. Consequently, the concept of efficiency becomes more relevant than optimality in multi-objective optimization problems. A feasible point $z^{*}$ is deemed to be an efficient solution of $(MOP_{BC})$ if there is no other $z \in \Feas$, such that both $\fvec(z) \leq \fvec(z^{*})$ and $\fvec(z) \neq \fvec(z^{*})$ hold, where $\fvec(z) = \left( \obj{1}(z), \obj{2}(z), \dots, \obj{m}(z) \right)$. Additionally, a feasible point $\xlm$ is characterized as a weak efficient solution of $(MOP_{BC})$ if there is no $z \in \Feas$ for which $\fvec(z) < \fvec(z^{*})$ is true. If $\Feas^*$ denotes the collection of efficient solutions of $(MOP_{BC})$, then $\fvec(\Feas^*)$ is termed the Pareto front of $(MOP_{BC})$.

\par If $z^*$ is a weak local efficient solution of $(MOP_{BC})$, then the following first-order necessary condition, which follows from Theorem 3.1.1 in \cite{Miettinen1998-to}, is satisfied at $z^*$. 
\begin{Theorem}[\textbf{First-order necessary condition for weak efficiency}]\label{theorem_2.1}
Let \(\xlm \in \Feas\) be a weak efficient solution of \((\mathrm{MOP}_{BC})\). Then there exist multipliers \((\mult, \mu) \in \Rplus{m} \times \Rplus{2n}\), \((\mult, \mu) \neq \mathbf{0}^{m+2n}\), such that
\begin{flalign}
\sum_{k = 1}^{m} \mult_k \nabla \obj{k}(\xlm) + \sum_{i = 1}^{n} \muc_i e_i - \sum_{i = n+1}^{2n} \muc_i e_i = 0;  \label{eq:FJ1} \\
\muc_i(\xlm_i-{ub}_i) = 0, \quad \forall~ i \in \Idx{n}; \label{eq:FJ2}\\
\muc_i ({lb}_i-\xlm_i) = 0, \quad \forall~ i \in \Idx{2n}\setminus \Idx{n}, \label{eq:FJ3}
\end{flalign}
where $e_1=(1,0,0,\dots,0)^T$. 
\end{Theorem}
\begin{definition}[\textbf{Critical point}]
A feasible point \(\xlm \in \Feas\) is said to be a \emph{critical point} of \((\mathrm{MOP}_{BC})\) if there exist multipliers \((\mult, \muc) \in \Rplus{m} \times \Rplus{2n}\) with \(\mult \neq \mathbf{0}^{m}\), which satisfy the conditions~\eqref{eq:FJ1}--\eqref{eq:FJ3}.
\end{definition}
\begin{note}
 If every objective function is convex, then a critical point is a weak efficient solution. A critical point is an efficient point if every objective function is a strictly convex function.
\end{note}
\section{Multi-objective tunneling method}\label{sec_tun}
In this section, we derive a multi-objective tunneling method for $(MOP_{BC})$ that can help us to find a global weak efficient solution for non-convex multi-objective optimization problems. The method proceeds through a sequence of cycles, each consisting of the following two phases:
\begin{enumerate}[label=(\Roman*)]
    \item the minimization phase, which aims to discover a critical point ($z^*$);
    \item the tunneling phase is designed to identify another critical point $\bar{z}\neq z^*$ and $\psi(\bar{z})<\psi(x^*)$ holds;
    \item the tunneling phase is executed repeatedly and stopped when\\ $\bar{z}=z^*$ holds and $z^*$ is considered as a global effecient solution of \((\mathrm{MOP}_{BC})\).

\end{enumerate}
Once finding out a local efficient solution, we construct a multi-objective variant of the auxiliary function, tunneling function as:
\begin{flalign*}
(\mathrm{TP}):~ \min_{z \in \mathbb{R}^n} \quad   &\left( \tunobj{1}(z), \tunobj{2}(z), \dots, \tunobj{m}(z) \right) \\
\text{s.t.} \quad & \tunobj{k}(z) \leq 0, \quad \forall~k \in \Idx{m
}, \\
& lb \leq z \leq ub,
\end{flalign*}
where,
\begin{equation*}
\tunobj{k}(z) =
\begin{cases} 
      \dfrac{\obj{k}(z) - \obj{k}(\xlm)}{\left[(z - \xlm)^\top (z - \xlm)\right]^{\eta}}, & \text{if } z \neq \xlm; \\
      \infty, & \text{otherwise};
\end{cases}
\end{equation*}
for all $k \in \Idx{m}$. Here, $\eta > 0$ is a prescribed constant and $\xlm$ denotes a weak efficient solution of $(\mathrm{MOP}_{BC})$.


One can observe that if $\bar{z}$ is a local efficient solution of $(\mathrm{TP})$, then $\obj{k}(\bar{z}) \leq \obj{k}(\xlm)$ for all $k \in \Idx{m}$ in a neighborhood of $\bar{z}$. The following theorem establishes that any critical point of $(\mathrm{TP})$ corresponds to a critical point of the original multi-objective problem $(\mathrm{MOP}_{BC})$.

\begin{Theorem}\label{theorem_3.1}
Let $\bar{z}$ be a critical point of the problem $(\mathrm{TP})$. 
If, for every direction $d$ and for step size $\alpha \in (0, \bar{\alpha}]$, 
the point $\bar{z} + \alpha d$ lies closer to $\xlm$ than $\bar{z}$, 
then $\bar{z}$ is also a critical point of the problem $(\mathrm{MOP}_{BC})$.
\end{Theorem}
\begin{proof}


Let $\bar{z}$ be a critical point for $(\mathrm{TP})$.  
Assume, for the sake of contradiction, that $\bar{z}$ is not a critical point for (MOP$_{BC}$).  
Then there exists a direction $d$ and a step size $\alpha \in (0, \bar{\alpha}]$ such that:
\[
\obj{j}(\bar{z} + \alpha d) < \obj{j}(\bar{z}) \quad \text{for all } j.
\]

Now consider the change in the $j$-th tunneling function:
\begin{align*}
\tunobj{j}(\bar{z} + \alpha d) - \tunobj{j}(\bar{z}) 
&= \frac{\obj{j}(\bar{z} + \alpha d) - \obj{j}(\xlm)}{\| \bar{z} + \alpha d - \xlm \|^{2\eta}} 
    - \frac{\obj{j}(\bar{z}) - \obj{j}(\xlm)}{\| \bar{z} - \xlm \|^{2\eta}} \\
&= \frac{\left[ \obj{j}(\bar{z} + \alpha d) - \obj{j}(\xlm) \right] \cdot \|\bar{z} - \xlm\|^{2\eta} 
     }%
    {\| \bar{z} + \alpha d - \xlm \|^{2\eta} \cdot \| \bar{z} - \xlm \|^{2\eta}}\\
    &\quad-\frac{\left[\obj{j}(\bar{z}) - \obj{j}(\xlm)\right] \cdot \| \bar{z} + \alpha d - \xlm \|^{2\eta}}{\| \bar{z} + \alpha d - \xlm \|^{2\eta} \cdot \| \bar{z} - \xlm \|^{2\eta}} \\
&= \frac{\obj{j}(\bar{z} + \alpha d) - \obj{j}(\bar{z})}{\| \bar{z} + \alpha d - \xlm \|^{2\eta}} \\
&\quad + \frac{\left[\obj{j}(\bar{z}) - \obj{j}(\xlm)\right] \cdot 
     \left[ \| \bar{z} - \xlm \|^{2\eta} - \| \bar{z} + \alpha d - \xlm \|^{2\eta} \right]}%
    {\| \bar{z} + \alpha d - \xlm \|^{2\eta} \cdot \| \bar{z} - \xlm \|^{2\eta}}.
\end{align*}




Since $\obj{j}(\bar{z} + \alpha d) < \obj{j}(\bar{z})$, the first term is strictly negative.  
Furthermore, because $\bar{z}$ is feasible in $(\mathrm{TP})$, we have \(\tunobj{j}(\bar{z}) \leq 0\), which implies:
$$\obj{j}(\bar{z}) - \obj{j}(\xlm) \leq 0.$$

In addition, the assumption that for every direction $d$ and step size $\alpha \in (0, \bar{\alpha}]$, 
the point $\bar{z} + \alpha d$ lies closer to $\xlm$ than $\bar{z}$ ensures:
\[
\| \bar{z} - \xlm \|^{2\eta} - \| \bar{z} + \alpha d - \xlm \|^{2\eta} \geq 0.
\]
Hence, the second term is non-positive.

Therefore, we conclude that:
\[
\tunobj{j}(\bar{z} + \alpha d) - \tunobj{j}(\bar{z}) < 0 \quad \text{for all } j,
\]
which implies that \(d\) is a descent direction of $\tun$ at $\bar{z}$, contradicting the assumption that \(\bar{z}\) is a critical point for the problem $(\mathrm{TP})$.

Hence, the assumption must be false. Therefore, $\bar{z}$ is also a critical point for \emph{(MOP\(_{BC}\))}.
\end{proof}
\begin{remark}
    From this theorem, one can observe that for $(\mathrm{MOP}_{BC})$, the tunneling phase makes it possible to obtain new critical points with improved objective values. Thus, by repeatedly applying the tunneling procedure, the method can escape from previously found critical points, which are not global efficient solutions and enrich the global Pareto front. This is justified in Section \ref{sec_exe} with different test problems.
\end{remark}
\subsection{Algorithm}
\label{sec_alg}
The above ideas are summarised in the following algorithm. In this algorithm, $PF$ denotes the Pareto front of $\mathrm(MOP_{BC})$ prior to tunneling, while $PFT$ signifies the Pareto front of $\mathrm(MOP_{BC})$ after tunneling. Furthermore, $WPF$ represents the weak Pareto front of $\mathrm(MOP_{BC})$ before tunneling and $WPFT$ refers to the weak Pareto front of $\mathrm(MOP_{BC})$ after tunneling.\\
\begin{algorithm}[H]
\caption{\textit{Tunneling method for nonlinear multi-objective optimization problems}}
\label{alg1}
\vspace{0.5mm}
\begin{algorithmic}
\State \textbf{Phase 0: Initialization}
\begin{enumerate}[label=(\alph*)]
  \item Supply $\fvec$ and a nonempty subset $\Feas^0\subset\Feas$ with $N$ initial points.
  \item Apply Steps~1--4 of Algorithm~6.1 in~\cite{ansary2021sqcqp} to modify these points and obtain $z^i$, $i=1,\dots,\ell$, satisfying inequality~(5.1) of the same reference, and compute ideal and nadir vectors via single objective global descent method~\cite{duanli2010global} for each $j\in\{1,\dots,m\}$.
  \item Set  $WPF=WPFT =PF =PFT = \emptyset$.
  \item Fix tolerance $\epsilon>0$ and an attempt cap $A_{\max}\in\mathbb{N}$.
\end{enumerate}
\State \textbf{for each} $z^i \in \{z^1,\dots,z^\ell\}$ \textbf{do} \;

\State \textbf{Phase 1: Minimization phase}
\begin{enumerate}[label=(\alph*)]
    \item Starting from $z^i$, apply some descent method to obtain a local weak efficient solution $z^{**}$ of $\fvec$ over $\Feas$; set $\xlm := z^{**}$.
  \item Update $WPF=WPF \cup \{\fvec(\xlm)\}$.
  \item Proceed to Phase 2.
\end{enumerate}
 \algstore{myalg}
\end{algorithmic}
\end{algorithm}
\begin{algorithm}[H]                     
\begin{algorithmic} [1]                   
\algrestore{myalg}
\State \textbf{Phase 2: Tunneling phase}
\begin{enumerate}[label=(\alph*)]
  \item Set anchor $z_{\mathrm{anc}} := \xlm$.
  \item \textbf{for} $a=1,2,\dots,A_{\max}$ \textbf{do}
    \begin{enumerate}[label=(b\arabic*)]
      \item Construct tunneling problem $(\mathrm{TP})$ with anchor $z_{\mathrm{anc}}$ and solve it to obtain a Fritz--John point $\bar{z}_1$.
      \item Construct $(\mathrm{TP})$ with anchor $\bar{z}_1$ and solve it to obtain a Fritz--John point $\bar{z}_2$.
      \item \textbf{If} $\|\fvec(\bar{z}_1)-\fvec(\bar{z}_2)\|<\epsilon$ \textbf{then} set $WPFT := WPFT \cup \{\fvec(\bar{z}_2)\}$ and \textbf{break}.
      \item \textbf{Else} set $z_{\mathrm{anc}} := \bar{z}_2$ and continue to the next attempt.
    \end{enumerate}
  \textbf{end for}
\end{enumerate}
\State \textbf{end for} \;
\State \textbf{Output:}
\begin{enumerate}[label=(\alph*)]
    \item Construct $PF$ by filtering nondominated points from $WPF$.
    \item Construct $PFT$ by filtering nondominated points from $WPFT$.
\end{enumerate}

\end{algorithmic}
\end{algorithm}

\section{Numerical illustration and discussion}
\label{sec_exe}

To evaluate the effectiveness of the proposed tunneling approach, we implemented a MATLAB (R2024b) script for the test problems mentioned in \linebreak
Table \ref{tab:table_1}. The following settings and procedures were used.
\begin{itemize}
\item \textit{Problem setup:} Consider a test function $\fvec(z)$ over the decision space \linebreak $[\text{lb}, \text{ub}]  \subset \Feas \subset \mathbb{R}^n$.    \item \textit{Finding ideal and nadir vector:} The ideal and nadir vectors are computed using the payoff table method \cite{Miettinen1998-to}. For each $\obj{j}$, global descent method for single objective optimization~\cite{duanli2010global} is applied using \texttt{fmincon} with three starting points: $\text{lb}$, $(\text{lb}+\text{ub})/2$ and $\text{ub}$.
    
    \item \textit{Generate initial points:}A total of 200 uniformly distributed random points are generated within the decision space \( [\text{lb}, \text{ub}] \subset \Feas \subset \mathbb{R}^n \) to initialize the optimization process.
    
    \item \textit{Spreading:} 
    Initial random points are refined using the spread subproblem $P_{\text{sp}}(z^0)$ as defined in section 5 of \cite{ansary2021sqcqp} solved by MATLAB function \texttt{quadprog}. Candidate points are retained if they satisfy inequality~(5.1) in \cite{ansary2021sqcqp}, ensuring well-distributed feasible solutions. This process continues until exactly $N=200$ spread-satisfying points are obtained.
        \item \textit{Minimization phase:} Each accepted spread point $z^i,i=1,2,\hdots,200$ is refined to find a local weak efficient solution by solving \texttt{MOSQCQP} method \cite{ansary2021sqcqp}. The resulting point $\xlm$ is added to the weak Pareto front $WPF$.

    \item \textit{Tunneling phase:}  
    At each refined solution $\xlm$, a tunneling problem is formulated as described in Section \ref{sec_tun}. The bounds are locally restricted to a neighborhood 
    \(\,[\max(\text{lb},~\xlm-\delta), \, \min(\text{ub},~\xlm+\delta)]\,\) with $\delta=10^{-5}$. Starting points are sampled randomly within this neighborhood, and \texttt{MOSQCQP} is used to minimize $\tun(z)$ to find $\bar{z}_1$, and again constructing the tunneling problem at $\bar{z}_1$ and repeating the process, another local weak efficient $\bar{z}_2$ solution is found.
    If $\|\fvec(\bar{z}_2)-\fvec(\bar{z}_1)\| < \epsilon$, where \linebreak
    $\epsilon=10^{-4}\sqrt{n}$, the process is deemed converged. Otherwise, the tunneling phase restarts up to maximum of six times with new random samples.

    \item \textit{Refinement:}
   After tunneling method, solutions are recorded in $WPFT$.   
    The final Pareto set $PF$ is obtained by removing the dominated solutions from $WPF$, while the refined Pareto set $PFT$ is constructed from $WPFT$ in a same manner.
\end{itemize}
\paragraph*{{\bf \em Experimental analysis}}
Algorithm~\ref{alg1} is tested on a collection of benchmark problems drawn from established references. See \cite{adhikary2025global,Evtushenko2013-ii,fliege2016sqp, hillermeier2001generalized,preuss2006pareto,schutze2008convergence, Zilinskas2014-po} for detailed descriptions. Each problem is initialized with 200 uniformly distributed random points and executed with the parameter choice $\eta = 0.9$. Table~\ref{tab:table_1} summarizes the computational setup and results, where $m$ and $n$ denote the number of objective functions and variables, respectively, while $\mathrm{MOSQCQP}$ and $\mathrm{MOTM}$ correspond to the number of approximate efficient solutions obtained using the sequential quadratic constrained quadratic programming method of \cite{ansary2021sqcqp} and the proposed multi-objective tunneling method.  

The practical performance of Algorithm~\ref{alg1} (MOTP) is evaluated on a broad set of box-constrained benchmark problems. The results indicate that the proposed approach is both robust and computationally efficient, producing consistent outcomes across diverse problem classes and, in many cases, outperforming established methods.  
\begin{table}[!htbp]
\centering
 \scalebox{0.79}{
\tiny  
\setlength{\tabcolsep}{6pt}  
\renewcommand{\arraystretch}{1.2}  
\begin{tabular}{|c|l|c|c|c|c|l|c|c|c|}

\hline
\textbf{Sl. no.} & \textbf{Test problem} & \textbf{(m,n)} & \textbf{MOSQCQP} & \textbf{MOTM} & 
\textbf{Sl. no.} & \textbf{Test problem} & \textbf{(m,n)} & \textbf{MOSQCQP} & \textbf{MOTM} \\
\hline
1	&	AL1	&	(2,20)	&	102	&	102	&	22	&	IKK1	&	(3,2)	&	200	&	200	\\
2	&	AL2	&	(2,50)	&	55	&	55	&	23	&	IM1	&	(2,2)	&	175	&	176	\\
3	&	CEC09\_1	&	(2,30)	&	188	&	191	&	24	&	Jin3	&	(2,2)	&	200	&	200	\\
4	&	CEC09\_2	&	(2,30)	&	198	&	200	&	25	&	Jin4\_a	&	(2,2)	&	173	&	174	\\
5	&	CEC09\_3	&	(2,30)	&	79	&	83	&	26	&	lovison2	&	(2,2)	&	177	&	180	\\
6	&	CEC09\_7	&	(2,30)	&	148	&	150	&	27	&	lovison3	&	(2,2)	&	161	&	161	\\
7	&	CEC09\_8	&	(3,30)	&	158	&	163	&	28	&	lovison4	&	(2,2)	&	185	&	185	\\
8	&	CL1	&	(2,2)	&	200	&	200	&	29	&	lovison5	&	(3,3)	&	84	&	84	\\
9	&	Deb521a\_a	&	(2,2)	&	200	&	200	&	30	&	lovison6	&	(3,3)	&	58	&	58	\\
10	&	Deb521b	&	(2,2)	&	200	&	200	&	31	&	MOP5	&	(3,2)	&	141	&	141	\\
11	&	DTLZ2	&	(3,12)	&	42	&	60	&	32	&	PNR1a	&	(2,2)	&	191	&	197	\\
12	&	DTLZ2n2	&	(2,2)	&	179	&	182	&	33	&	PNR1b	&	(2,2)	&	200	&	200	\\
13	&	DTLZ5	&	(3,12)	&	111	&	112	&	34	&	PNR1c	&	(2,2)	&	178	&	179	\\
14	&	DTLZ5n2	&	(2,2)	&	196	&	198	&	35	&	PNR1d	&	(2,2)	&	200	&	200	\\
15	&	DTLZ6	&	(3,22)	&	157	&	158	&	36	&	PNR1e	&	(2,2)	&	200	&	200	\\
16	&	DTLZ6n2	&	(2,2)	&	195	&	193	&	37	&	Shekel	&	(2,2)	&	182	&	182	\\
17	&	EP2	&	(2,2)	&	200	&	200	&	38	&	slcdt1	&	(2,2)	&	200	&	200	\\
18	&	EX005	&	(2,2)	&	185	&	199	&	39	&	VFM1	&	(3,2)	&	178	&	200	\\
19	&	LP1	&	(2,50)	&	43	&	43	&	40	&	ZDT1	&	(2,30)	&	170	&	171	\\
20	&	GE5	&	(3,3)	&	200	&	200	&	41	&	ZDT2	&	(2,30)	&	200	&	200	\\
21	&	hil	&	(2,2)	&	158	&	158	&	42	&	ZDT3	&	(2,30)	&	115	&	115	\\\hline
\end{tabular}
 }
\caption{Comparison of (\textit{MOTM}) and (\textit{MOSQCQP}) for various test problems.}
\label{tab:table_1}
\end{table}
Table~\ref{tab:table_1} provides a detailed comparison between the MOSQCQP method and the MOTM across 42 benchmark problems of varying dimensions. The results show that MOTM is broadly competitive and frequently superior. In 23 problems, both algorithms yield the same number of nondominated solutions, underscoring their comparable reliability. In 18 cases, however, MOTM achieves clear improvements by recovering a larger portion of the Pareto front, while MOSQCQP surpasses MOTM in only one instance, and by a negligible margin.  

Across the benchmark suite, MOTM consistently enhances Pareto coverage, with the scale of improvement depending on problem complexity. On simpler test cases, it provides modest yet reliable gains, whereas in more complex nonlinear problems, the improvements are much more pronounced. Representative examples such as \texttt{DTLZ2}, \texttt{EX005}, and \texttt{VFM1} highlight MOTM’s ability to uncover broader Pareto fronts, demonstrating its effectiveness in avoiding premature convergence and capturing a richer spectrum of trade-off solutions. MOTM also exhibits superior robustness by attaining the maximum coverage of $200$ nondominated points in 14 benchmarks, compared with 12 for MOSQCQP. Such consistent recoveries across problem classes, ranging from convex structures to irregular landscapes, highlight its adaptability and stability. Overall, the findings confirm that MOTM not only matches MOSQCQP on a large fraction of problems but also delivers consistent and often substantial improvements elsewhere. Its ability to achieve broader and denser Pareto approximations, together with a higher rate of full-front recoveries, establishes MOTM as a reliable and effective framework for practical multi-objective optimization.
\paragraph*{{\bf \em Performance analysis}}
In this part, Algorithm~\ref{alg1} (MOTM) is benchmarked against the modified SQCQP method proposed in~\cite{ansary2021sqcqp} using \textit{performance profiles} with respect to different metrics. For more details on performance profiles, the readers may refer to \cite{ansary2020sequential, fliege2016sqp, more2009benchmarking, zitzler2003performance}. \\
\textit{Performance profiles:} A performance profile is a cumulative distribution function, denoted as $\rho(\tau)$, that captures the relative efficiency of solvers over a set of test problems, with respect to a specified performance metric.\\

Formally, consider a set of solvers $\mathcal{SO}$ and a collection of test problems $\mathcal{P}$. Let $\varsigma_{p,s}$ denote the performance of solver $s$ on problem $p$. The performance ratio is defined as $r_{p,s} = \varsigma_{p,s} / \min_{s \in \mathcal{SO}} \varsigma_{p,s}$, representing the efficiency of solver $s$ relative to the best solver on problem $p$. The cumulative performance profile for solver $s$ is then defined as:
\begin{equation*}
 \rho_s(\tau) = \frac{|\{p \in \mathcal{P} : r_{p,s} \leq \tau\}|}{|\mathcal{P}|}
\end{equation*}
This represents the fraction of test problems on which a given solver achieves results within a factor $\tau$ of the best known performance. In multi-objective optimization, where algorithms generate a set of non-dominated solutions rather than single optima, standard comparison metrics are often insufficient. To effectively assess the quality and distribution of the obtained Pareto fronts, additional performance metrics are incorporated into the construction of the performance profile.\\
\textit{Purity metric:} Let $F_{p,s}$ denote the approximated Pareto front for problem $p$ obtained using method $s$. An approximation to the true Pareto front $F_p$ can be constructed by taking the union of all approximated fronts across the methods, $\bigcup_{s \in S} F_{p,s}$ and removing all dominated solutions. To evaluate the contribution of each method to this reference set, the \textit{purity metric} is employed. For a given algorithm $s$ and problem $p$, the purity is defined as:
\begin{equation*}
\tilde{t}_{p,s} := \frac{|F_{p}|}{|F_p \cap F_{p,s}|}.
\end{equation*}
This metric reflects the inverse proportion of non-dominated solutions produced by method $s$ relative to the total number of non-dominated solutions in the reference front. $\tilde{t}_{p,s} = \infty$ implies that the method failed to generate any non-dominated solution in the reference Pareto front for problem $p$.\\
\textit{Spread metrics:} Two types of spread metrics, denoted by $\Gamma$ and $\Delta$, are utilized to evaluate whether the solutions generated by a given method are well-distributed within the approximated Pareto front of a particular problem. Let $x^1, x^2, \ldots, x^N$ denote the set of solutions obtained by method $s$ for problem $p$ and let these be ordered such that $f_j(x^i) \leq f_j(x^{i+1})$ for $i = 1, 2, \ldots, N - 1$. Let $x^0$ and $x^{N+1}$ denote, respectively, the best-known approximations to the global minimum and global maximum of $f_j$, computed over all approximated Pareto fronts across methods. Define $\delta_{i,j} = f_j(x^{i+1}) - f_j(x^i)$. The $\Gamma$-spread metric is then given by:
\begin{equation*}
   \Gamma_{p,s} := \max_{j \in \Lambda_m} ~\max_{i \in \{0,1,\ldots,N\}} \delta_{i,j}.
\end{equation*}
Let $\bar{\delta}_j$ denote the average of $\delta_{i,j}$ for $i = 1, 2, \ldots, N - 1$. The $\Delta$ spread metric for method $s$ on problem $p$ is defined as
\begin{equation*}
    \Delta_{p,s} := \max_{j \in \Lambda_m} \left( \frac{\delta_{0,j} + \delta_{N,j} + \sum_{i=1}^{N-1} |\delta_{i,j} - \bar{\delta}_j|}{\delta_{0,j} + \delta_{N,j} + (N-1)\bar{\delta}_j} \right).
\end{equation*}
\textit{Hypervolume metric:} The hypervolume metric quantifies the extent of the objective space dominated by a given approximation of the Pareto front with respect to a predefined reference point \( P_{\text{ref}} \). It serves as a comprehensive quality measure by capturing both convergence and diversity of the solution set. However, due to the computational complexity of determining the exact hypervolume, particularly in higher dimensions, various approximation strategies have been proposed in the literature. In the present study, the hypervolume is estimated through Monte Carlo sampling, wherein 10{,}000 points are uniformly generated within the hyper-rectangle defined by the reference point \( P_{\text{ref}} \) and the ideal point of the objective space. The metric value is then computed as $hv_{p,s} = \frac{N_{\text{dom}}}{10{,}000}$, where \( N_{\text{dom}} \) denotes the number of sampled points that are dominated by the approximate Pareto front. A higher value of \( hv_{p,s} \) reflects superior performance using hypervolume metric. For performance profiling and solver comparisons, the inverse form \( \widetilde{hv_{p,s}} = 1 / hv_{p,s} \) is employed to ensure consistency with minimization-based metric frameworks.

\textit{Function evaluations:} 
To quantify the computational effort of each method, we assess performance using the number of function evaluations. For each test problem, the union of all Pareto fronts obtained across different methods is taken, and the reference front is determined by removing all dominated points. The contribution of each method is then measured by computing the intersection of its approximated Pareto front with the reference front.  

Let $F_{p,s}$ denote the Pareto front obtained by method $s$ on problem $p$, and $F_p$ the reference front. We define the function evaluation metric as:
\[
   f_{\text{eval}}(p,s) := \frac{\text{total number of function evaluations by } s}{|F_{p,s} \cap F_p|}.
\]
If $|F_{p,s} \cap F_p| = 0$, we assign $f_{\text{eval}}(p,s) = \infty$, indicating that method $s$ failed to contribute to the reference front.
\par The above metric is then incorporated into performance profiles to enable a fair comparison of computational efficiency across all methods.

\begin{figure}[H]
    \centering
    \begin{subfigure}[b]{0.45\textwidth}
        \centering
        \includegraphics[width=\textwidth,height=2.5cm]{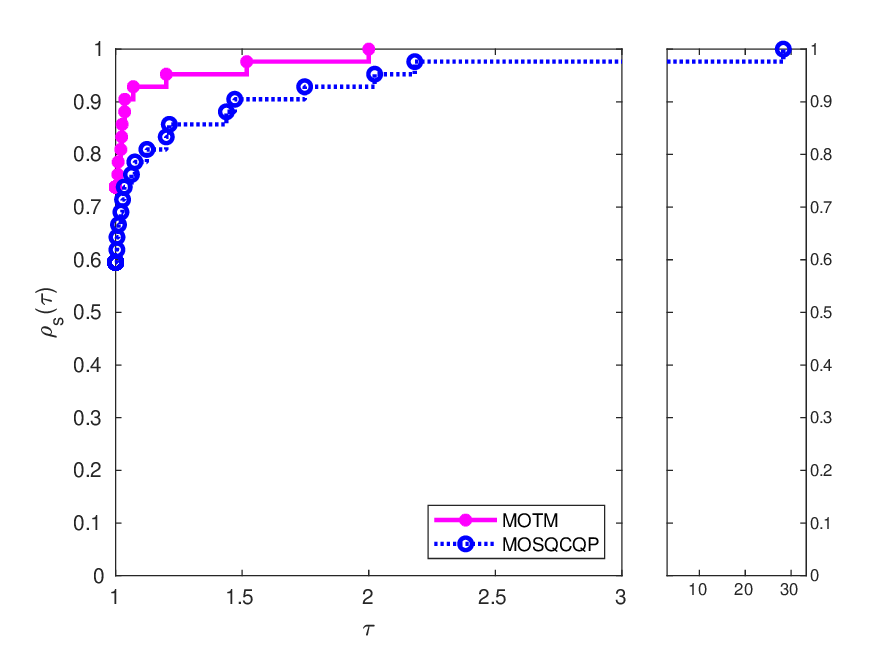}
        \caption{\raggedright Performance profile between MOTM and MOSQCQP}
        \label{fig_purity_motp_mosqcqp}
    \end{subfigure}
    \hfill
    \begin{subfigure}[b]{0.45\textwidth}
        \centering
        \includegraphics[width=\textwidth,height=2.5cm]{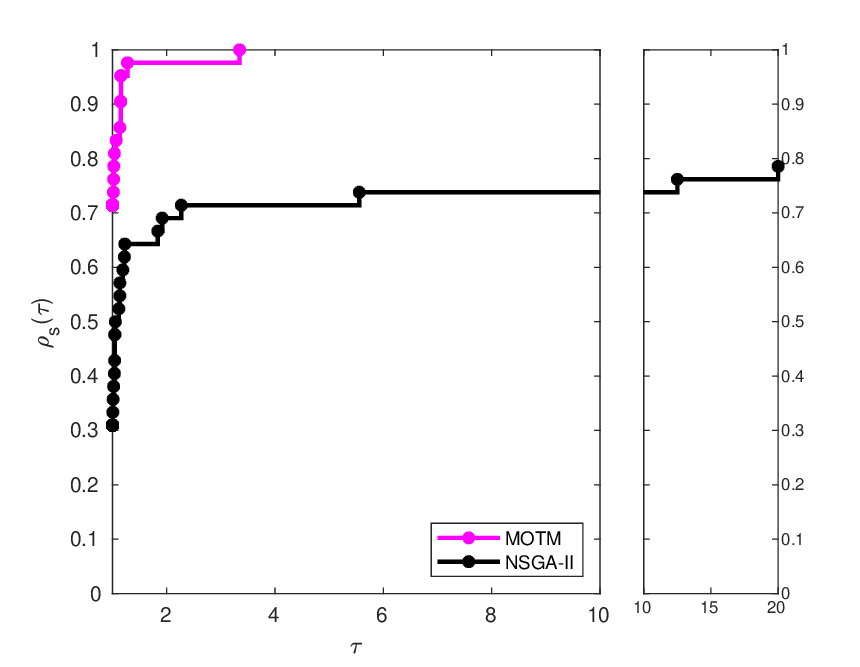}
        \caption{\raggedright Performance profile between MOTM and NSGA-II}
        \label{fig_purity_motp_nsga2}
    \end{subfigure}    
    \caption{Performance profiles using purity metric}
    \label{fig_purity}
\end{figure}
\begin{figure}[!http]
    \centering
    \begin{subfigure}[b]{0.45\textwidth}
        \centering
        \includegraphics[width=\textwidth,height=2.5cm]{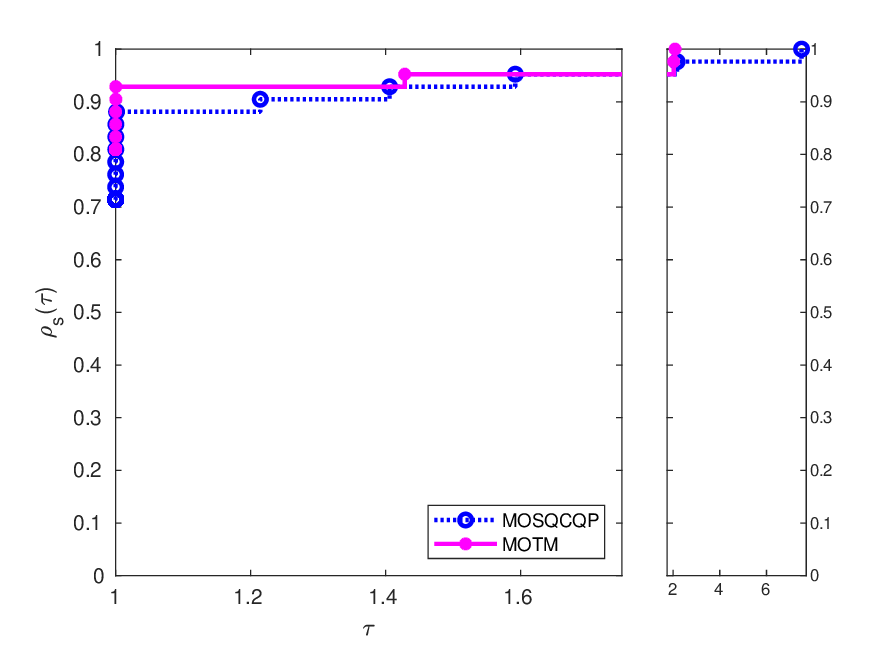}
        \caption{\raggedright Performance profile between MOTM and MOSQCQP}
        \label{fig_gamma_motp_mosqcqp}
    \end{subfigure}
    \hfill
    \begin{subfigure}[b]{0.45\textwidth}
        \centering
        \includegraphics[width=\textwidth,height=2.5cm]{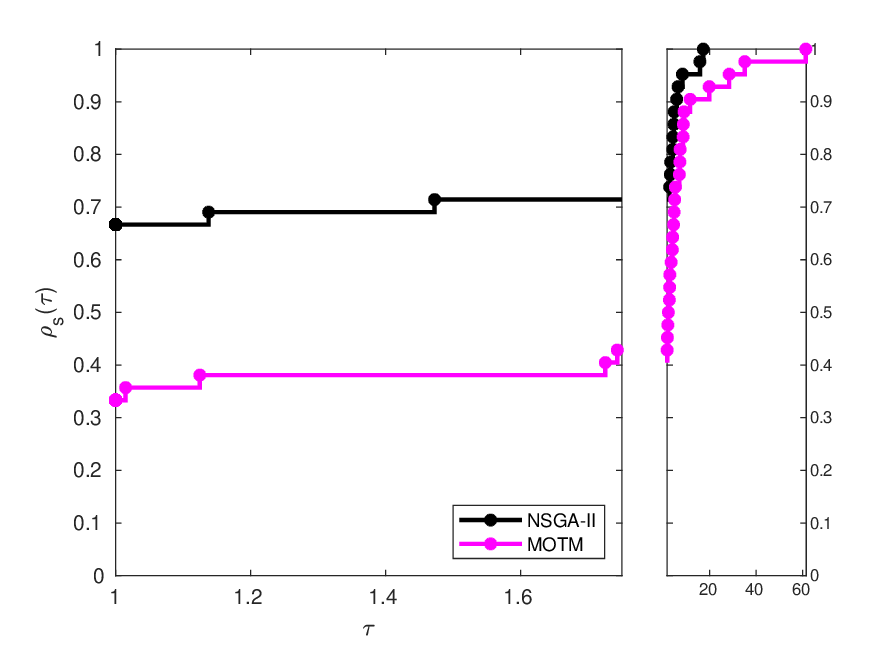}
        \caption{\raggedright Performance profile between MOTM and NSGA-II}
        \label{fig_gamma_motp_nsga2}
    \end{subfigure}    
    \caption{Performance profiles using $\Gamma$-spread metric}
    \label{fig_gamma}
\end{figure}

\begin{figure}[H]
    \centering
    \begin{subfigure}[b]{0.45\textwidth}
        \centering
        \includegraphics[width=\textwidth,height=2.5cm]{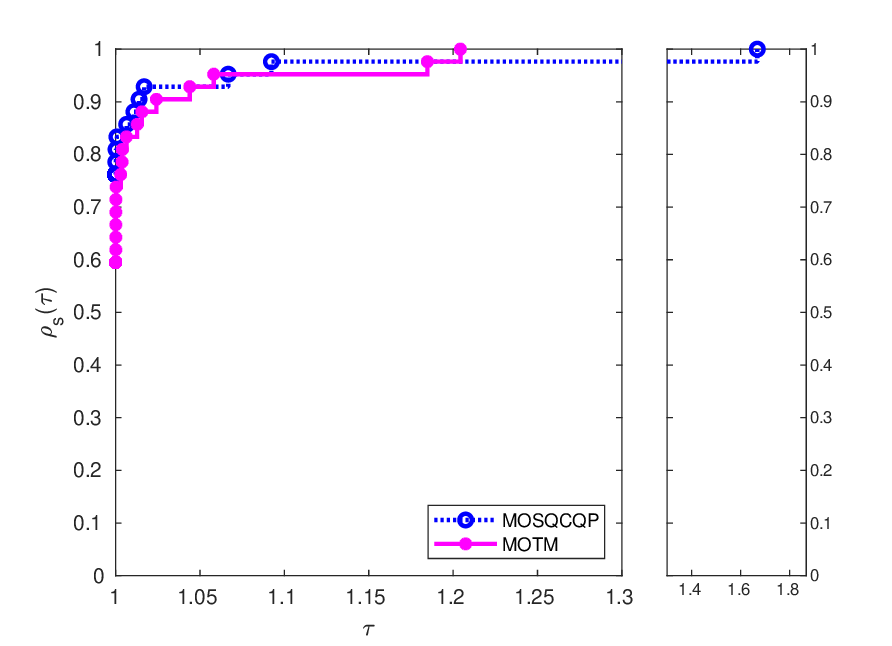}
        \caption{\raggedright Performance profile between MOTM and MOSQCQP}
        \label{fig_delta_motp_mosqcp}
    \end{subfigure}
    \hfill
    \begin{subfigure}[b]{0.45\textwidth}
        \centering
        \includegraphics[width=\textwidth,height=2.5cm]{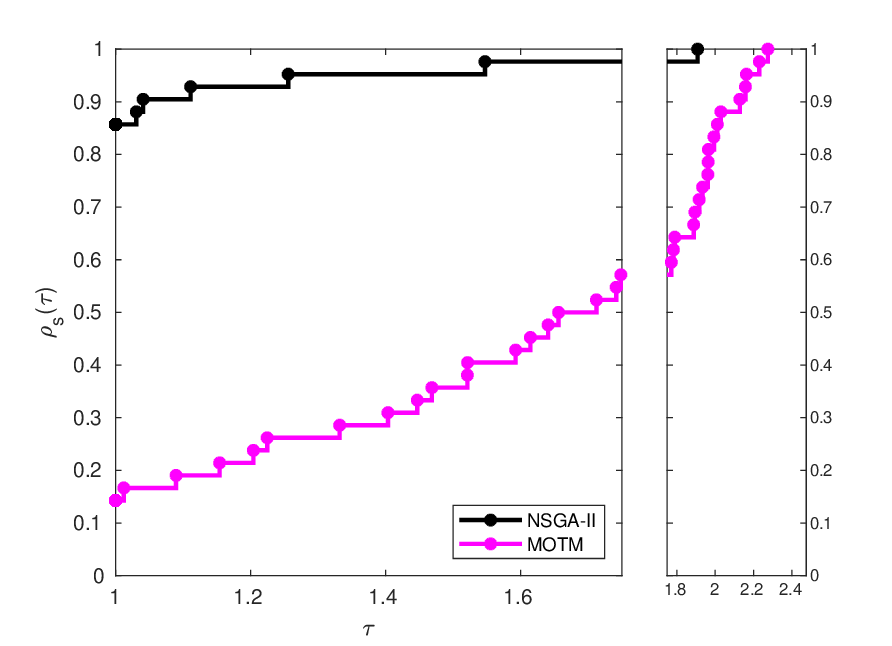}
        \caption{\raggedright Performance profile between MOTM and NSGA-II}
        \label{fig_delta_motp_nsga2}
    \end{subfigure}    
    \caption{Performance profiles using $\Delta$-spread metric}
    \label{fig_delta}
\end{figure}


\begin{figure}[H]
    \centering
    \begin{subfigure}[b]{0.45\textwidth}
        \centering
        \includegraphics[width=\textwidth,height=2.5cm]{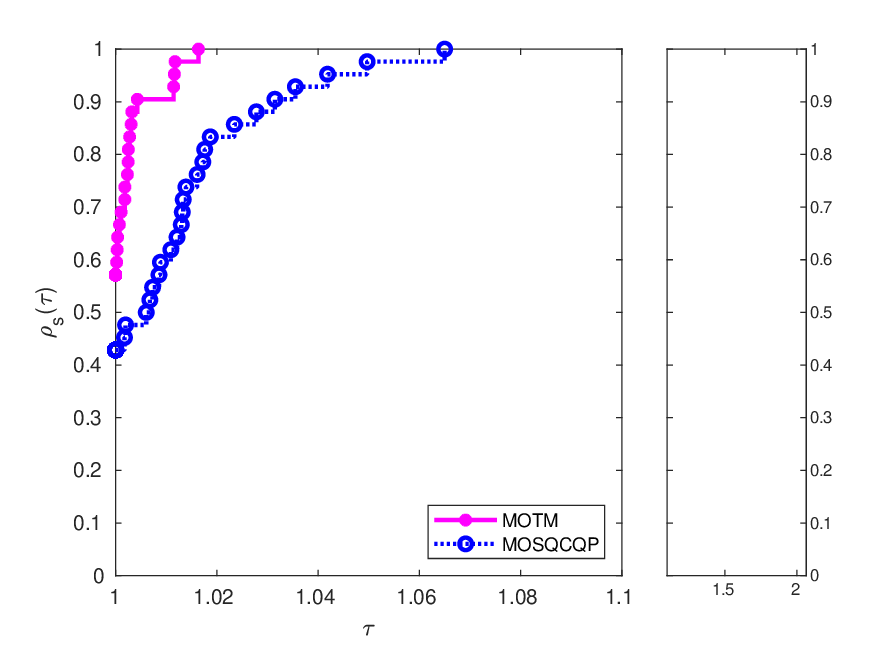}
        \caption{\raggedright Performance profile between MOTM and MOSQCQP}
        \label{fig_hv_motp_mosqcqp}
    \end{subfigure}
    \hfill
    \begin{subfigure}[b]{0.45\textwidth}
        \centering
        \includegraphics[width=\textwidth,height=2.5cm]{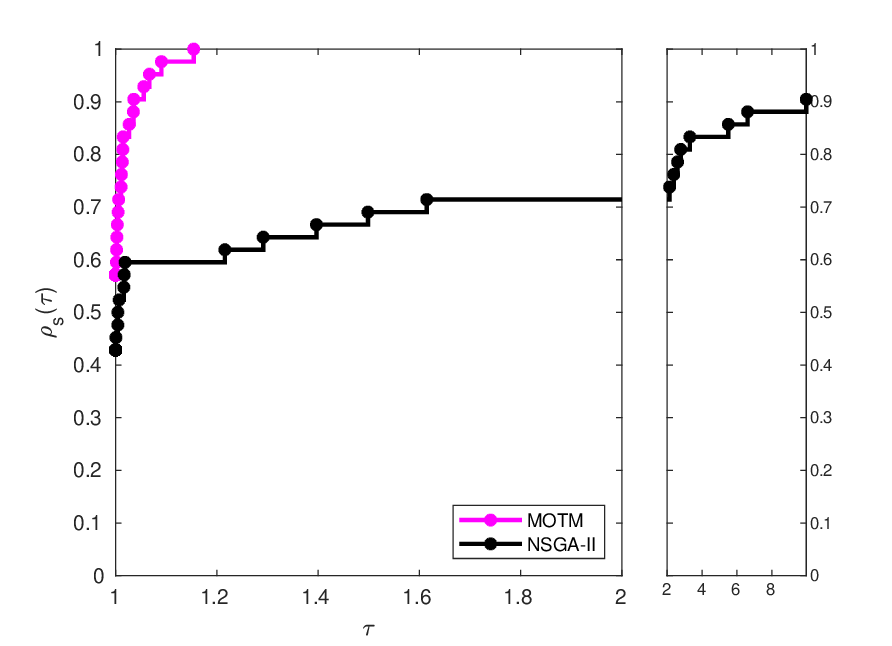}
        \caption{\raggedright Performance profile between MOTM and NSGA-II}
        \label{fig_hv_motp_nsga2}
    \end{subfigure}    
    \caption{Performance profiles using hypervolume metric}
    \label{fig_hv}
\end{figure}


\begin{figure}[H]
    \centering
    \begin{subfigure}[b]{0.45\textwidth}
        \centering
        \includegraphics[width=\textwidth,height=2.5cm]{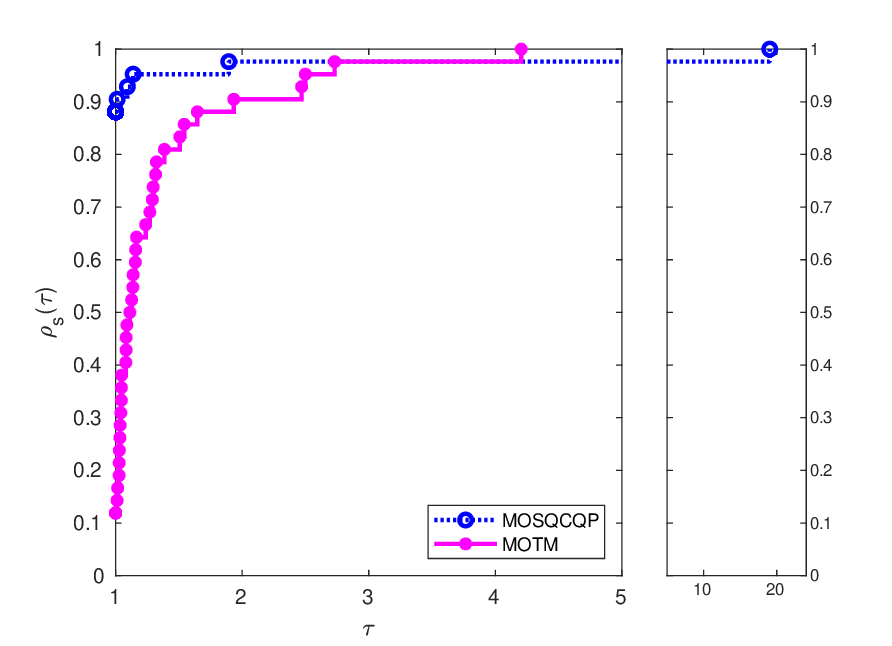}
        \caption{\raggedright Performance profile between MOTM and MOSQCQP}
        \label{fig_feval_motp_mosqcqp}
    \end{subfigure}
    \hfill
    \begin{subfigure}[b]{0.45\textwidth}
        \centering
        \includegraphics[width=\textwidth,height=2.5cm]{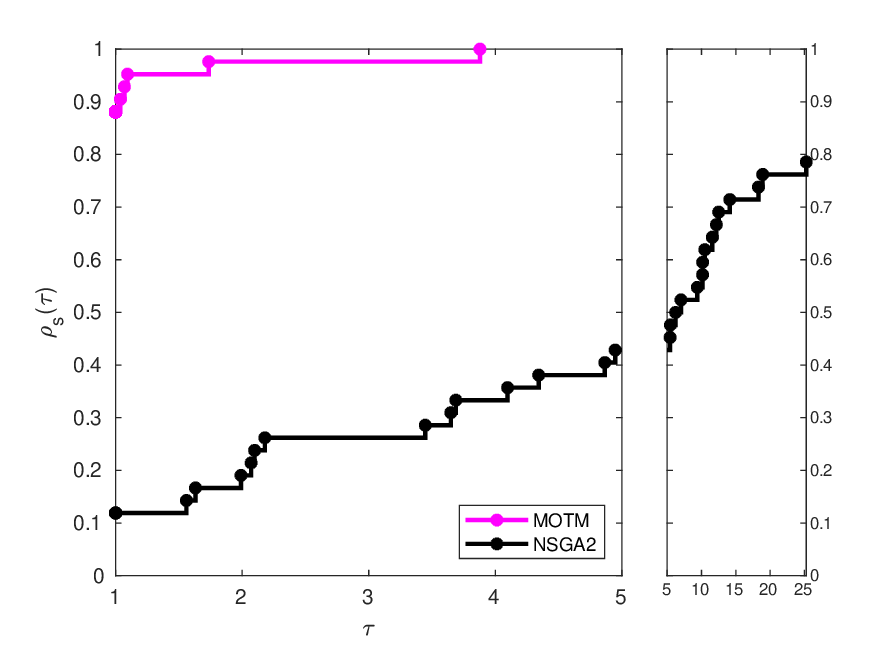}
        \caption{\raggedright Performance profile between MOTM and NSGA-II}
        \label{fig_feval_motp_nsga2}
    \end{subfigure}    
    \caption{Performance profiles using number of function evaluations}
    \label{fig_feval}
\end{figure}
To benchmark the efficiency and robustness of the proposed multi-objective tunneling method (MOTM), we compare it with two established baselines: the sequential quadratic constrained quadratic programming framework \linebreak
(MOSQCQP)~\cite{ansary2021sqcqp} and the evolutionary algorithm NSGA-II~\cite{deb2002fast}. The evaluation is conducted across the benchmark problems reported in Table~\ref{tab:table_1}, with performance profiles constructed using five complementary indicators: the purity index, the $\Delta$-spread metric, the $\Gamma$-spacing measure, the hypervolume index, and the number of function evaluations. This suite of metrics jointly captures convergence accuracy, distribution quality, solution uniqueness, and computational cost~\cite{ansary2020sequential,ansary2021sqcqp,fliege2016sqp}.  

The purity profiles reveal that MOTM contributes a larger proportion of unique nondominated solutions compared to both MOSQCQP and NSGA-II. This indicates that MOTM consistently uncovers novel Pareto-optimal points that competing methods miss, thereby enriching the global approximation set.  

The $\Delta$-spread results (Figure~\ref{fig_delta}) present a more nuanced outcome. Figure~\ref{fig_delta_motp_mosqcp} shows that MOSQCQP holds a slight advantage under strict thresholds, yielding marginally more uniform spacing. However, this advantage diminishes with relaxed tolerances, where MOTM achieves comparable spread quality. In contrast, Figure~\ref{fig_delta_motp_nsga2} highlights that NSGA-II produces superior spread, with a higher and faster-rising profile reflecting more evenly distributed solutions.  

The $\Gamma$-spreading profiles complement these findings by assessing consistency of solution spacing. Here, MOTM performs competitively with \linebreak
MOSQCQP but trails NSGA-II, whose distribution is more uniform across the front. While this suggests scope for improving MOTM’s spacing regularity, its other strengths largely offset this limitation.  

The hypervolume profiles (Figure~\ref{fig_hv}) consistently favor MOTM, which achieves broader Pareto front coverage than both MOSQCQP and NSGA-II. This reflects stronger convergence to the true front and an enhanced ability to explore diverse trade-off regions.  

Finally, the efficiency profiles (Figure~\ref{fig_feval}) illustrate the computational cost. Since MOSQCQP serves as an intermediate process, it requires fewer evaluations, whereas the tunneling-based exploration in MOTM incurs a higher expense. Nonetheless, MOTM remains competitive and consistently outperforms NSGA-II, delivering superior Pareto approximations at a lower cost than evolutionary search.  

Taken together, these results highlight complementary strengths of MOTM. Although somewhat more expensive than MOSQCQP, it provides broader Pareto coverage, higher purity, and superior hypervolume performance, while remaining competitive in distribution quality. Against NSGA-II, MOTM consistently delivers better convergence, coverage, and efficiency, even though it trails in strict spread and spacing metrics. Overall, this balance of robustness, convergence, and computational efficiency establishes MOTM as a reliable and competitive framework for solving challenging multi-objective optimization problems.
\section{Conclusion}\label{(sec_con)}
In this paper, we have developed a multi-objective tunneling method for nonlinear box constrained multi-objective optimization problems. This method is free from any kind of priori chosen parameter or ordering information of objective functions. The tunneling algorithm showcases promising results in enriching the solution space and getting better efficient solutions for various nonlinear multi-objective optimization problems. While the current work is limited to smooth unconstrained problems, future research will focus on extending the framework to constrained and nonsmooth problems, thereby enhancing its scope and applicability to complex multi-objective optimization problems.
\section*{Acknowledgement} The authors thank the referees for their detailed comments and suggestions that have significantly improved the content as well as the presentation of the results in the paper.
\section*{Declarations}
\begin{itemize}
\item {\bf Funding:} This article is partially funded by {\it prime minister's research fellowship (PMRF)} scheme. Author Bikram Adhikary acknowledges funding from the PMRF scheme (PMRF ID. 2202747), India.
\item {\bf Competing interests:} The authors declare no competing interests.
\item {\bf Ethics approval and consent to participate:} Not applicable.
\item {\bf Consent for publication:} Not applicable.
\item {\bf Data availability:} In this paper, we have not used any associated data. 
\item {\bf Materials availability:} The simulation codes are available upon reasonable request from the corresponding author.
\item {\bf Author contribution:} Both authors have contributed equally to theoretical developments and numerical calculations.
\end{itemize}
 \bibliographystyle{imsart-number}

\end{document}